\documentclass[12pt]{amsart}
\usepackage{amssymb}
\usepackage[all]{xy}

\newcommand{\grbl}{{\mbox{\textit{\tiny gp}}}}
\long\def\forget#1\forgotten{}
\newcommand{\issuenumber}{18}
\newcommand{\issuemonth}{September}
\newcommand{\issueyear}{2006}

\newtheorem{thm}{Theorem}[section]
\newtheorem{prob}[thm]{Problem}

\newtheorem{issue}{Issue}

\theoremstyle{definition}

\theoremstyle{remark}

\newcommand{\ed}{
\general\end{document}}

\newcommand{\Cantor}{{\{0,1\}^\N}}

\newcommand{\fb}{\mathfrak{b}}

\newcommand{\fd}{\mathfrak{d}}
\newcommand{\fp}{\mathfrak{p}}

\newcommand{\NON}{{\mathsf   {NON}}}
\newcommand{\COF}{{\mathsf   {COF}}}

\newcommand{\M}{\mathcal{M}}

\newcommand{\cov}{\mathsf{cov}}
\newcommand{\add}{\mathsf{add}}
\newcommand{\cof}{\mathsf{cof}}

\newcommand{\non}{\mathsf{non}}

\newcommand{\CH}{the Continuum Hypothesis}
\newcommand{\R}{\mathbb{R}}

\newcommand{\fo}{\mathfrak{od}}

\renewcommand{\split}{\mathsf{Split}}
\newcommand{\bq}{\begin{quote}}
\newcommand{\eq}{\end{quote}}
\newcommand{\cO}{\mathcal{O}}
\newcommand{\B}{\mathcal{B}}
\newcommand{\BG}{\B_\Gamma}

\newcommand{\BO}{\B_\Omega}

\newcommand{\sone}{\mathsf{S}_1}    \newcommand{\sfin}{\mathsf{S}_{fin}}

\newcommand{\ufin}{\mathsf{U}_{fin}}

\newcommand{\nin}{\not\in}

\newcommand{\N}{\mathbb{N}}

\newcommand{\sbst}{\subseteq}
\newcommand{\by}[2]{\par\hfill\emph{#1}, #2}
\newcommand{\nby}[1]{\par\hfill\emph{#1}}
\newcommand{\Tau}{\mathrm{T}}
\newcommand{\CE}{\textsc{CE}}

\newcommand{\be}{\begin{enumerate}}
\newcommand{\ee}{\end{enumerate}}
\newcommand{\bi}{\begin{itemize}}
\newcommand{\ei}{\end{itemize}}
\renewcommand{\i}{\item}

\newcommand{\general}{\small\vfill\par\noindent\hrulefill\par
\noindent\textbf{Previous issues.} The previous issues of this
bulletin, which contain general information (first issue), basic
definitions, research announcements, and open problems (all
issues) are available online,
at \texttt{http://front.math.ucdavis.edu/search?\&t=\%22SPM+Bulletin\%22}
\\[0.1cm]
\textbf{Contributions.}
Please submit your contributions (announcements, discussions, and open problems)
by e-mailing us. It is preferred to write them
in \LaTeX{}.
The authors are urged to use as standard notation as possible, or otherwise give
the definitions or a reference to where the notation is explained.
Contributions to this bulletin would not require any transfer of copyright,
and material presented here can be published elsewhere.\\[0.1cm]
\textbf{Subscription.}
To receive this bulletin (free) to your
e-mailbox, e-mail us:\\
{boaz.tsaban@weizmann.ac.il}
}

\newcommand{\nArxPaper}[5]{\subsection{#2}{#4}\par\hfill{\arx{#1}}\par\hfill\emph{#3}}

\newcommand{\nAMSPaper}[5]{\subsection{#2}{#4}\par\hfill{\texttt{#1}}\par\hfill\emph{#3}}

\newcommand{\arx}[1]{\texttt{http://arxiv.org/math/#1}}
\newcommand{\url}[1]{\bq\texttt{#1}\eq}
\newcommand{\online}[1]{The paper is available online at \url{#1}}

\title[$\mathcal{SPM}$ Bulletin \textbf{\issuenumber} (\issuemonth{} \issueyear)]{%
$\mathcal{SPM}$ Bulletin\\[0.5cm]
Issue number \issuenumber: \issuemonth{} \issueyear{} \CE{}}

\begin{document}
\maketitle

\tableofcontents

\section{Editor's note}

Good news in brief:
\be
\i The abstracts for the recent \emph{Kielce Conference on Set-Theoretic Topology}
are available online at
\url{http://www.pu.kielce.pl/~topoconf/skrzynia/abstrakty.pdf}
and at
\url{http://atlas-conferences.com/cgi-bin/abstract/catc-01}
Some of the presentations of this beautiful conference should be
available soon in
\url{http://www.pu.kielce.pl/~topoconf/prezentacje.html}
\i Some pictures from the recent TOPOSYM and Kielce conferences
are now available in Marcin Kysiak's webpage
\url{http://www.mimuw.edu.pl/\~{}mkysiak}
(follow the links at the bottom of this page).
\i July 14-19, 2007: The \emph{Logic Colloquium 2007} (ASL European Summer Meeting)
will take place in Wroc\l{}aw (Poland), one of the world's centers for SPM related
mathematics. It is likely that some of the talks there will be related to SPM.
\ee

\medskip

Contributions to the next issue are, as always, welcome.

\medskip

\by{Boaz Tsaban}{boaz.tsaban@weizmann.ac.il}

\hfill \texttt{http://www.cs.biu.ac.il/\~{}tsaban}

\section{Research announcements}

\nAMSPaper{http://www.ams.org/journal-getitem?pii=S0002-9939-06-08371-7}
{A surprising covering of the real line}
{G\'abor Kun}
{We construct an increasing sequence of Borel subsets of $\R$, such
that their union is $\R$, but $\R$ cannot be covered with countably many translations
of one set. The proof uses a random method.}

\nAMSPaper{http://dx.doi.org/10.1016/S0166-8641(01)00106-7}
{Unions of chains in dyadic compact spaces and topological groups}
{Mikhail G.\ Tkachenko and Yolanda Torres Falc\'on}
{The following problem is considered: If a topological group $G$ is the union of an increasing chain
of subspaces and certain cardinal invariants of the subspaces are known, what can be said about $G$?
We prove that if $G$ is locally compact and every subspace in the chain has countable pseudocharacter
or tightness, then $G$ is metrizable. We also prove a similar assertion for $\sigma$-compact and totally
bounded groups represented as the union of first countable subspaces, when the length of the chain
is a regular cardinal greater than $\aleph_1$.
Finally, we show that these results are not valid in general, not
even for compact spaces.}

\nArxPaper{0606270}
{On the Pytkeev property in spaces of continuous functions}
{Petr Simon and Boaz Tsaban}
{Answering a question of Sakai, we show that the minimal cardinality of a set of reals $X$ such that $C_p(X)$
does not have the Pytkeev property is equal to the pseudo-intersection number $\fp$.
Our approach leads to a natural characterization of the Pytkeev property of $C_p(X)$ by means of a covering
property of $X$, and to a similar result for the Reznicenko property of $C_p(X)$.
We also give a new result of Miller: If $C_p(X)$ has the Pytkeev property, then $X$ has strong measure zero.
(To appear in \emph{Proceedings of the AMS}.)}

\nArxPaper{math.GN/0608107}
{Selection principles related to $\alpha_i$-properties}
{Ljubi\v{s}a D.R.\ Ko\v{c}inac}
{We investigate selection principles which are motivated by
Arhangel'ski\v{i}'s $\alpha_i$-properties, $i=1,2,3,4$, and their
relations with classical selection principles. It will be shown
that they are closely related to the selection principle $\sone$
and often are equivalent to it.}

\nArxPaper{0607592}
{On the Kocinac $\alpha_i$ properties}
{Boaz Tsaban}
{The Kocinac $\alpha_i$ properties, $i=1,2,3,4$, are generalizations of Arkhangel'skii's
$\alpha_i$ local properties. We give a complete classification of these properties when
applied to the standard families of open covers of topological spaces or to the standard
families of open covers of topological groups.
One of the latter properties characterizes totally bounded groups.
We also answer a question of Kocinac.}

\nArxPaper{0606285}
{A new selection principle}
{Boaz Tsaban}
{Motivated by a recent result of Sakai, we define a new selection operator for
covers of topological spaces, inducing new selection hypotheses. We initiate a
systematic study of the new hypotheses. Some intriguing problems remain open.}

\nArxPaper{math.GN/0608035}
{First Countable Continua and Proper Forcing}
{Joan E. Hart and Kenneth Kunen}
{Assuming the Continuum Hypothesis, there is a compact first countable connected space
of weight $\aleph_1$ with no totally disconnected perfect subsets. Each such space, however,
may be destroyed by some proper forcing order which does not add reals.}

\nArxPaper{math.GN/0608086}
{The convergence space of minimal usco mappings}
{R Anguelov, O.\ F.K.\ Kalenda}
{A convergence structure generalizing the order convergence structure on the set of
Hausdorff continuous interval functions is defined on the set of minimal usco maps.
The properties of the obtained convergence space are investigated and essential links
with the pointwise convergence and the order convergence are revealed.
The convergence structure can be extended to a uniform convergence structure so that the
convergence space is complete. The important issue of the denseness of the subset of all
continuous functions is also addressed.}

\nArxPaper{math.GN/0609090}
{$D$-forced spaces: a new approach to resolvability}
{Istvan Juhasz, Lajos Soukup, and Zoltan Szentmiklossy}
{We introduce a ZFC method that enables us to build spaces (in
fact special dense subspaces of certain Cantor cubes) in which we
have ``full control'' over all dense subsets. Using this method we
are able to construct, in ZFC, for each uncountable regular
cardinal $\lambda$ a 0-dimensional $T_2$, hence Tychonov, space
which is $\mu$-resolvable for all $\mu < \lambda$ but not
$\lambda$-resolvable. This yields the final (negative) solution of
a celebrated problem of Ceder and Pearson raised in 1967: Are
$\omega$-resolvable spaces maximally resolvable? This method
enables us to solve several other open problems concerning
resolvability as well.

The paper appeared in Top. Appl. \textbf{153} (2006), 1800--1824.}

\newcommand{\ps}{\operatorname{ps}}
\newcommand{\pe}{\operatorname{pe}}
\nArxPaper{math.GN/0609091}
{Resolvability of spaces having small spread or extent}
{Istvan Juhasz, Lajos Soukup, and Zoltan Szentmiklossy}
{In a recent paper O. Pavlov proved the following two interesting resolvability results:
\be
\i If a space $X$ satisfies $\Delta(X) > \ps(X)$ then $X$ is maximally resolvable.
\i If a $T_3$-space $X$ satisfies $\Delta(X) > \pe(X)$ then $X$ is $\omega$-resolvable.
\ee

Here $\ps(X)$ ($\pe(X)$) denotes the smallest successor cardinal
such that $X$ has no discrete (closed discrete) subset of that
size and $\Delta(X)$ is the smallest cardinality of a non-empty
open set in $X$. In this note we improve (1) by showing that
$\Delta(X) >$ $\ps(X)$ can be relaxed to $\Delta(X) \ge$ $\ps(X)$.
In particular, if $X$ is a space of countable spread with
$\Delta(X) > \omega$ then $X$ is maximally resolvable. The
question if an analogous improvement of (2) is valid remains open,
but we present a proof of (2) that is simpler than Pavlov's.}

\nArxPaper{math.GN/0609092}
{Resolvability and monotone normality}
{Istvan Juhasz, Lajos Soukup, and Zoltan Szentmiklossy}
{A space $X$ is said to be $\kappa$-resolvable (resp. almost
$\kappa$-resolvable) if it contains $\kappa$ dense sets that are
pairwise disjoint (resp. almost disjoint over the ideal of nowhere
dense subsets). $X$ is maximally resolvable iff it is
$\Delta(X)$-resolvable, where $\Delta(X) = \min\{|G| : G \ne
\emptyset {open}\}.$
We show that every crowded monotonically
normal (in short: MN) space is $\omega$-resolvable and almost
$\mu$-resolvable, where $\mu = \min\{2^{\omega}, \omega_2 \}$. On
the other hand, if $\kappa$ is a measurable cardinal then there is
a MN space $X$ with $\Delta(X) = \kappa$ such that no subspace of
$X$ is $\omega_1$-resolvable. Any MN space of cardinality $<
\aleph_\omega$ is maximally resolvable. But from a supercompact
cardinal we obtain the consistency of the existence of a MN space
$X$ with $|X| = \Delta(X) = \aleph_{\omega}$ such that no subspace
of $X$ is $\omega_2$-resolvable.}

\nAMSPaper{http://www.ams.org/journal-getitem?pii=S0002-9939-06-08542-X}
{Isomorphism of Borel full groups}
{Benjamin D.\ Miller; Christian Rosendal}
{Suppose that $G$ and $H$ are Polish groups which act in a Borel
fashion on Polish spaces $X$ and $Y$. Let $E_G^X$ and $E_H^Y$
denote the corresponding orbit equivalence relations, and $[G]$
and $[H]$ the corresponding Borel full groups. Modulo the obvious
counterexamples, we show that $[G] \cong [H]$ iff $E_G^X$ is Borel
isomorphic to $E_H^Y$.}

\newcommand{\qkappa}{{\mathbb Q}(\kappa)}
\nArxPaper{math.LO/0608642}
{A Poset Hierarchy}
{M. D{\v{z}}amonja and K. Thompson}
{This article extends a paper of Abraham and Bonnet which
generalised the famous Hausdorff characterisation of the class of
scattered linear orders. Abraham and Bonnet gave a poset hierarchy
that characterised the class of scattered posets which do not have
infinite antichains (abbreviated FAC for finite antichain
condition). An antichain here is taken in the sense of
incomparability. We define a larger poset hierarchy than that of
Abraham and Bonnet, to include a broader class of ``scattered''
posets that we call $\kappa$-scattered. These posets cannot embed
any order such that for every two subsets of size $ < \kappa$, one
being strictly less than the other, there is an element in
between. If a linear order has this property and has size $\kappa$
we call this set $\qkappa$. Such a set only exists when
$\kappa^{<\kappa}=\kappa$. Partial orders with the property that
for every $a<b$ the set $\{x: a<x<b\}$ has size $\geq \kappa$ are
called weakly $\kappa$-dense, and partial orders that do not have
a weakly $\kappa$-dense subset are called strongly
$\kappa$-scattered. We prove that our hierarchy includes all
strongly $\kappa$-scattered FAC posets, and that the hierarchy is
included in the class of all FAC $\kappa$-scattered posets. In
addition, we prove that our hierarchy is in fact the closure of
the class of all $\kappa$-well-founded linear orders under
inversions, lexicographic sums and FAC weakenings. For
$\kappa=\aleph_0$ our hierarchy agrees with the one from the
Abraham-Bonnet theorem.

Central European Journal of Mathematics, vol. 4, no. 2, (2006), 225--241.
}

\nArxPaper{math.FA/0608616}
{Infinite asymptotic games}
{Christian Rosendal}
{We study infinite asymptotic games in Banach spaces with an
F.D.D. and prove that analytic games are determined by
characterising precisely the conditions for the players to have
winning strategies. These results are applied to characterise
spaces embeddable into $\ell_p$ sums of finite dimensional spaces,
extending results of Odell and Schlumprecht, and to study various
notions of homogeneity of bases and Banach spaces. These results
are related to questions of rapidity of subsequence extraction
from normalised weakly null sequences.}

\nArxPaper{math.GN/0608376}
{Elementary submodels and separable monotonically normal compacta}
{Todd Eisworth}
{In this note, we use elementary submodels to prove that a
separable monotonically normal compactum can be mapped on a
separable metric space via a continuous function whose fibers have
cardinality at most $2$.}

\nArxPaper{math.LO/0608382}
{An application of CAT}
{Mirna D\v{z}amonja and Jean Larson}
{We comment on a question of Justin Moore on colorings of pairs
of nodes in an Aronszajn tree and solve an instance of it.}

\nArxPaper{math.LO/0608384}
{A general Stone representation theorem}
{Mirna D\v{z}amonja}
{This note contains a Stone-style representation theorem for compact Hausdorff spaces.}

\newcommand{\algB}{\mathfrak B}
\nArxPaper{math.LO/0608336}
{Measure Recognition Problem}
{Mirna D\v{z}amonja}
{This is an article on the example of the Measure Recognition Problem (MRP).
The article highlights the phenomenon of the utility of a multidisciplinary
mathematical approach to a single mathematical problem, in
particular the value of a set-theoretic analysis. MRP asks if for
a given Boolean algebra $\algB$ and a property $\Phi$ of measures
one can recognize by purely combinatorial means if $\algB$
supports a strictly positive measure with property $\Phi$. The
most famous instance of this problem is MRP(countable additivity),
and in the first part of the article we survey the known results
on this and some other problems. We show how these results
naturally lead to asking about two other specific instances of the
problem MRP, namely MRP(nonatomic) and MRP(separable). Then we
show how our recent work D\v zamonja and Plebanek (2006) gives an
easy solution to the former of these problems, and gives some
partial information about the latter. The long term goal of this
line of research is to obtain a structure theory of Boolean
algebras that support a finitely additive strictly positive
measure, along the lines of Maharam theorem which gives such a
structure theorem for measure algebras.}

\nArxPaper{math.GN/0609351}
{Cardinal invariants for $C$-cross topologies}
{Andrzej Kucharski and Szymon Plewik}
{$C$-cross topologies are introduced. Modifications of the Kuratowski-Ulam Theorem are considered.
Cardinal invariants $\add$, $\cof$, $\cov$, and $\non$ with respect to meager or nowhere dense
subsets are compared. Remarks on invariants $\cof(nwd_Y)$ are mentioned for dense subspaces $Y$ of $X$.}

\section{Problem of the Issue}

In our work (in progress) on Ramsey theory of open covers,
we have encountered the following natural problem, to which
we do not see the answer.

\begin{prob}\label{prm}
Assume that $X\sbst\R$ is Hurewicz and all finite powers of $X$ are Menger.
Are all finite powers of $X$ Hurewicz?
\end{prob}

In modern notation \cite{coc7}, Problem \ref{prm} can be stated as follows.

\begin{prob}\label{2}
Is $\ufin(\cO,\Gamma)\cap \sfin(\Omega,\Omega)=\sfin(\Omega,\Omega^\grbl)$?
\end{prob}

$\ufin(\cO,\Gamma)=\sfin(\Omega,\Lambda^\grbl)$ \cite{coc7}, which in
turn is equivalent to $\binom{\Lambda}{\Lambda^\grbl}$ \cite{hureslaloms}
(see also \cite{SF1}).

Thus, another way to state this is:

\begin{prob}\label{3}
Is it true that, for $X\sbst\R$ satisfying $\sfin(\Omega,\Omega)$,
$\binom{\Lambda}{\Lambda^\grbl}=\binom{\Omega}{\Omega^\grbl}$?
\end{prob}

A negative answer under \CH{} is what we would like to have, so here
is a variant:

\begin{prob}[CH]\label{ch}
Is there a Hurewicz $X$ such that $X^2$ is Menger but not Hurewicz?
\end{prob}

\nby{Nadav Samet and Boaz Tsaban}

\subsection{Solution to Problem \ref{prm}} 
Problem \ref{prm} (and thus also Problems \ref{2} and \ref{3})
was solved in the negative, by Scheepers and Tall, in their paper
\emph{Lindel\"of indestructibility, topological games and selection principles},
Fundamenta Mathematicae 210 (2010), 1--46.

Problem \ref{ch} remains open.

\newpage

\section{Unsolved problems from earlier issues}

\begin{issue}
Is $\binom{\Omega}{\Gamma}=\binom{\Omega}{\Tau}$?
\end{issue}

\begin{issue}
Is $\ufin(\Gamma,\Omega)=\sfin(\Gamma,\Omega)$?
And if not, does $\ufin(\Gamma,\Gamma)$ imply
$\sfin(\Gamma,\Omega)$?
\end{issue}

\stepcounter{issue}

\begin{issue}
Does $\sone(\Omega,\Tau)$ imply $\ufin(\Gamma,\Gamma)$?
\end{issue}

\begin{issue}
Is $\fp=\fp^*$? (See the definition of $\fp^*$ in that issue.)
\end{issue}

\begin{issue}
Does there exist (in ZFC) an uncountable set satisfying $\sone(\BG,\B)$?
\end{issue}

\stepcounter{issue}

\begin{issue}
Does $X \nin \NON(\M)$ and $Y\nin\mathsf{D}$ imply that
$X\cup Y\nin \COF(\M)$?
\end{issue}

\begin{issue}[CH]
Is $\split(\Lambda,\Lambda)$ preserved under finite unions?
\end{issue}

\begin{issue}
Is $\cov(\M)=\fo$? (See the definition of $\fo$ in that issue.)
\end{issue}

\begin{issue}
Does $\sone(\Gamma,\Gamma)$ always contain an element of cardinality $\fb$?
\end{issue}

\begin{issue}
Could there be a Baire metric space $M$ of weight $\aleph_1$ and a partition
$\mathcal{U}$ of $M$ into $\aleph_1$ meager sets where for each ${\mathcal U}'\subset\mathcal U$,
$\bigcup {\mathcal U}'$ has the Baire property in $M$?
\end{issue}

\stepcounter{issue} 

\begin{issue}
Does there exist (in ZFC) a set of reals $X$ of cardinality $\fd$ such that all
finite powers of $X$ have Menger's property $\sfin(\cO,\cO)$?
\end{issue}

\begin{issue}
Can a Borel non-$\sigma$-compact group be generated by a Hurewicz subspace?
\end{issue}

\begin{issue}[MA]
Is there an uncountable $X\sbst\R$ satisfying $\sone(\BO,\BG)$?
\end{issue}

\begin{issue}[CH]
Is there a totally imperfect $X$ satisfying $\ufin(\cO,\Gamma)$
that can be mapped continuously onto $\Cantor$?
\end{issue}

\begin{issue}[CH]
Is there a Hurewicz $X$ such that $X^2$ is Menger but not Hurewicz?
\end{issue}

\ed
\begin{thebibliography}{00}

\bibitem{coc7}
Lj.\ D.R.\  Ko\v{c}inac and M.\ Scheepers,
\emph{Combinatorics of open covers (VII): Groupability},
Fundamenta Mathematicae \textbf{179} (2003), 131--155.

\bibitem{hureslaloms}
B.\ Tsaban,
\emph{The Hurewicz covering property and slaloms in the Baire space},
Fundamenta Mathematicae \textbf{181} (2004), 273--280.

\bibitem{SF1}
L.\ Zdomskyy,
\emph{A semifilter approach to selection principles},
Commentationes Mathematicae Universitatis Carolinae \textbf{46} (2005), 525--540.

\end{thebibliography}
